\documentstyle[12pt]{article}
\newcommand{\scr}[1]{\mbox{\scriptsize $#1$}}

\newcommand{\sm}[1]{\mbox{\small $#1$}}
\newcommand{\la}[1]{\mbox{\large $#1$}}
\newcommand{\La}[1]{\mbox{\Large $#1$}}
\newcommand{\LA}[1]{\mbox{\LARGE $#1$}}

\newcommand{\lnab}[1]{\la{\nabla}_{\!\!#1}}

\newcommand{\lnabo}[1]{\la{\nabla}^{\bot}_{\!\!#1}} 

\newcommand{\qed}{\mbox{~~~\boldmath $\Box$}} 
\newcommand{\R}[1]{I\!\!R^{#1}}
\newcommand{\Co}{l\!\!\! C}
\newcommand{\ra}{\rightarrow}
\newcommand{\non}{\nonumber}
\newcommand{\ha}{\scr{\frac{1}{2}}}

\newcommand{\al}{\alpha}
\newcommand{\be}{\beta}

\newcommand{\bbe}{\bar{\beta}}

\newcommand{\bmu}{\bar{\mu}}
\newcommand{\brho}{\bar{\rho}}

\newcommand{\gdf}[3]{g(\mbox{\large ${\nabla}$}_{\!\!{#1}} dF({#2}),JdF({#3}))}

\newcommand{\Jw}{J_{\omega}}
\newcommand{\kw}{{\cal K}_{\omega}}

\newcommand{\Fw}{F^{*}\omega}
\newcommand{\wt}{\omega^{\bot}}
\newcommand{\Jt}{J^{\bot}}
\newtheorem{Lm}{Lemma}
\newtheorem{Pp}{Proposition}
\newtheorem{Th}{Theorem}
\newtheorem{Cr}{Corollary}

\begin{document}
\baselineskip .6cm
\setcounter{page}{1}
\title{ On the K\"{a}hler angles of Submanifolds\\
{\small To the memory of Giorgio Valli}}
 \author{
Isabel M.\ C.\ Salavessa}
\date{}
\maketitle ~~~\\[-5mm]
\mbox{ } \\[-10mm]
{\small {\bf Abstract:}
 We prove that under certain conditions on the mean curvature and on the
  K\"{a}hler angles, a compact submanifold $M$
of real dimension $2n$,  immersed into a K\"{a}hler-Einstein  manifold $N$ of complex dimension $2n$,  must be either a 
complex or a Lagrangian submanifold of $N$, or have constant K\"{a}hler
angle, depending on $n=1$, $n=2$, or $n\geq 3$, and the sign of the scalar 
curvature of $N$. These results generalize to non-minimal submanifolds
some known results for minimal submanifolds.
Our main tool is a Bochner-type tecnique involving a
formula on the Laplacian of a symmetric function on the K\"{a}hler angles
and the Weitzenb\"{o}ck formula for the K\"{a}hler form of $N$ restricted to $M$.
}\\[1mm]
\markright{\sl\hfill  Salavessa \hfill} 
\section{Introduction}
\setcounter{Th}{0}
\setcounter{Pp}{0}
\setcounter{Cr} {0}
\setcounter{Lm} {0}
\setcounter{Def} {0}
\setcounter{equation} {0}
Let $(N,J,g)$ be a K\"{a}hler-Einstein manifold of complex dimension $2n$, 
complex structure $J$,  Riemannian metric $g$, and 
 $F:M^{2n}\ra N^{2n}$ be 
an immersed submanifold $M$ of real dimension $2n$.
We denote by $\omega(X,Y)=g(JX,Y)$ the  K\"{a}hler form
and by $R$ the scalar curvature of $N$, that is, the Ricci tensor
of $N$ is given by $Ricci = Rg$. The cosine of the  K\"{a}hler angles $\{\theta_{\al}\}_{1\leq \al\leq n}$ are the eigenvalues of  $\Fw$.
If  the eigenvalues  are all equal to 0 (resp. 1), $F$
is a  Lagrangian (resp. complex) submanifold.  A natural question is to ask
if $N$ allows submanifolds with arbitrary given K\"{a}hler angles
and  mean curvature. 
An answer is that, the K\"{a}hler angles and  the second fundamental 
form of $F$, and the Ricci tensor of $N$ are interrelated. 
Conditions on some of these
geometric objects have implications for the other ones. 
There are obstructions to the existence of minimal Lagrangian submanifolds 
in a general K\"{a}hler manifold, but these obstructions do not occur in
a  K\"{a}hler\--Einstein manifolds, where such submanifolds exist with abundance 
([Br]). This is the reason we choose K\"{a}hler-Einstein manifolds as
ambient spaces.
An example how the sign of the scalar curvature of $N$ determines the 
K\"{a}hler angles is the fact that
if $F$ is a totally geodesic immersion  and $N$ is not Ricci-flat, then either 
$F$ has a complex direction, or $F$ is Lagrangian ([S-V,1]).
A relation among the $\theta_{\al}$, $\lnab{}dF$, and $R$ can be described 
through a formula on the Laplacian of a locally Lipschitz map $\kappa$,
symmetric on the K\"{a}hler angles of $F$,  where the Ricci tensor of $N$
and some components of the second fundamental form of $F$ appear.
Such  kind of formula was used for minimal immersions
in [W,1] for $n=1$, and in [S-V,1,2] for $n\geq 2$.  

 A natural condition
for $n\geq 2$  is
to impose equality on the K\"{a}hler angles.
Products of surfaces immersed with the same constant K\"{a}hler angle
$\theta$ into
 K\"{a}hler-Einstein surfaces of the same scalar curvature $R$, give
  submanifolds immersed with constant equal K\"{a}hler angle $\theta$
into a K\"{a}hler-Einstein manifold of scalar curvature $R$.
The slant submanifolds introduced and exhaustively studied by B-Y Chen 
(see e.g. [Che,1,2], [Che-M], [Che-T,1,2]) are just submanifolds with constant and equal K\"{a}hler
angles. Examples are given in complex spaces form,  some of them via Hopf's
fibration [Che-T,1,2].
A minimal 4-dimensional submanifold of a Calabi-Yau manifold of
complex  dimension 4,  calibrated by a Cayley calibration,  also called
Cayley submanifold, 
 is just the same
as a  minimal submanifold with
 equal K\"{a}hler angles ([G]). Existence theory  of such submanifolds in $\Co^{4}$, with given boundary data, is guaranteed by the  theory of calibrations of Harvey and Lawson [H-L]. 

Submanifolds with equal K\"{a}hler angles have a  role in 4 and 8 dimensional gauge
theories. For example,  each of such Cayley submanifolds in
$\Co^{4}$ carries a 21-dimensional family of (anti)-self-dual $SU(2)$ Yang-Mills fields [H-L]. Recentely, Tian [T] proved that blow-up loci
of complex anti-self-dual instantons on  Calabi-Yau 4-folds are
Cayley cycles, which are, except for a set of 4-dimensional
Hausdorff measure  zero, a countable
union of $C^1$  4-dimensional Cayley submanifolds.

If $N$ is an hyper--K\"{a}hler manifold of
complex dimension 4 and hyper-K\"{a}hler structure $(J_{x})_{x\in S^2}$,
any  submanifold of real dimension 4 that is $J_x$-complex for some
$x\in S^2$, is a minimal submanifold  with equal K\"{a}hler angles
of each $(N,J_{y},g)$ ([S-V,2]), and the common K\"{a}hler angle is given by $\cos\theta(p)=\|(J_{y}X)^{\top}\|$,
where $X$ is any unit vector of $T_{p}M$. A proof of this assertion is simply 
to remark that, if $\{X, J_xX, Y, J_x Y\}$ is an o.n. basis of $T_pM$,
then the matrix of the K\"{a}hler form $\omega_y$ w.r.t.\ $J_y$,  restricted
to this basis, is just a multiple of a matrix in $\R{4}$ that represents
an orthogonal complex structure of $\R{4}$, i.e.\  of the type
$aI+bJ+cK$, where $I,J,K$ defines the usual hyper-K\"{a}hler structure
of $\R{4}$, and $a^2+b^2+c^2=1$. The square of this multiple is given by
$\langle x,y\rangle^2+\langle J_y X,Y\rangle^2 +\langle J_{x\times y}X,Y\rangle^2=
\|(J_y X)^{\top}\|^2$.
 This example suggests us
a way to build examples of (local) submanifolds with
equal K\"{a}hler angles. Let $(N,I,g)$ be a K\"{a}hler manifold of complex dimension 4, and $U\subset N$  an open set where an orthornormal
frame  of the form $\{X_1,IX_1, X_2, IX_2, Y_1,IY_1, Y_2, IY_2 \}$ is
defined. If for each $p\in U$, we identify $T_p N$ with $\R{4}\times \R{4}$, 
through this frame,  we are defining a family 
of local $g$-orthogonal almost complex structures $ J_{x}= a i\times i+b j\times j+c k\times k$, for $x=(a,b,c) \in S^2$, where
$i,j,k$ denotes de canonical hyper-K\"{a}hler structure of $\R{4}$. 
Then any almost $J_x$-complex 4-dimensional submanifold  $M$ 
is a submanifold with equal K\"{a}hler angles of the K\"{a}hler manifold 
$(N, I, g)$. It may not be minimal, because $J_x$ may
not be a K\"{a}hler structure, or not even integrable. 

Such a condition on the K\"{a}hler angles, turns out to
be more restrictive for sub\-ma\-ni\-folds of non Ricci-flat manifolds, or if 
$M$ is closed, that is, compact and orientable. 
A  combination of the formula of $\triangle\kappa$
for minimal immersions with equal K\"{a}hler angles, 
 with the Weitzenb\"{o}ck formula for
$\Fw$, lead us in [S-V,2]  to the  conclusion that the K\"{a}hler angle
must be constant, and in general it is either $0$ or $\frac{\pi}{2}$.
Namely, we have:
\begin{Th}  Let $F:M^{2n}\ra N^{2n}$ be a minimal immersion with  equal K\"{a}hler angles.\\[1mm]
$(i)~(${\bf [W,1]}$)$ If $n=1$, $M$ is closed, 
$R<0$, and $F$ has no
complex points, then $F$ is Lagrangian.\\[1mm]
$(ii)~(${\bf [S-V,2], [G]}$)$ If $n=2$ and $R\neq 0$, then  
$F$ is either a complex 
or a Lagrangian submanifold.\\[1mm]
$(iii)~(${\bf [S-V,2]}$)$ If $n\geq 3$,  $M$ is closed, and   
$R<0$,  then $F$ is either a complex or a Lagrangian submanifold.\\
$(iv)~(${\bf [S-V,2]}$)$ If $n\geq 3$,  $M$ is closed, $R=0$,  
then the  common K\"{a}hler angle  must be constant. 
\end{Th}
If $n=2$ and $R=0$ we cannot conclude the K\"{a}hler angle is constant.
It is easy to find examples of
minimal immersions with constant and non-constant 
equal K\"{a}hler angle, for the case of $M$ not compact and $N$ 
the Euclidean space. Namely, the most simple family of submanifolds with constant equal K\"{a}hler angle
of $~\Co^{2n}$ can be given by the vector subspaces defined by a linear map
$F:\R{2n}\ra \Co^{2n}\equiv (\R{2n}\times \R{2n}, J_{0})$,
$F(X)=(X, a \Jw X)$, where $a$ is any  real number and $\Jw$
is a $g_0$-orthogonal complex structure of $\R{2n}$, and where $g_0$ is 
the Euclidean metric and $J_0(X,Y)=(-Y,X)$. These are totally geodesic submanifolds
with constant equal K\"{a}hler angle $\cos\theta=\frac{2|a|}{1+a^2}$, and
$\Fw(X,Y)= \cos\theta \,F^{*}\!g_{0}\,(\pm\Jw X,Y)$, with $F^{*}\!g_{0}$ a 
$\Jw$-hermitian euclidean metric.
In  ([D-S]) we have the following example of non-constant K\"{a}hler angle
well away from $0$. The graph of  the anti-$i$-holomorphic map $f:\R{4}\ra \R{4}$  given by $f(x,y,z,w)=(u,v,-u,-v)$, where
\[\begin{array}{l}
u(x,y,z,w)=\phi(x+z)\xi '(y+w),\\
v(x,y,z,w)=-\phi'(x+z)\xi(y+w)\\
\phi(t)=\sin t,~~~\xi(t)=\sinh t, \end{array}\]
defines a minimal complete submanifold of $\Co^{4}$ with equal K\"{a}hler
angles satisfying
\[\cos\theta =\frac{2\sqrt{\cos^2(x+z)+\sinh^2(y+w)}}{1+4(\cos^2(x+z)+\sinh^2(y+w))}.\] 
This graph has no complex points, for
$0\leq \cos\theta\leq \ha$, and the set of Lagrangian points  is a 
infinite discrete union of disjoint 2-planes, 
\[{\cal L}=\bigcup_{-\infty\leq k\leq +\infty} span_{\R{}}\{(1,0,-1,0), (0,1,0,-1)\} +(0,0,(\ha+k)\pi,0).\]

In this paper  we present a formula for $\triangle\kappa$,
but now not
assuming minimality of $F$, obtaining some extra terms involving the
mean curvature $H$ of $F$. 
We will see that the  above conclusions still hold for $F$ not minimal, but
under certain weaker condition on the mean curvature of $F$. These conclusions
show how rigid  K\"{a}hler-Einstein manifolds are with respect to the K\"{a}hler angles and the mean curvature of a submanifold, leading to some non-existence  of certain types of submanifolds, depending on the sign of
the scalar curvature $R$ of $N$ and on the dimension $n$.

We summarize the main results of this paper:
\begin{Th} Assume  $n=2$, and  $M$ is closed,
$N$ is non Ricci-flat,  and $F:M\ra N$
is an immersion with equal K\"{a}hler angles, 
$\theta_{\al}=\theta ~\forall\al$. If  
\begin{equation}
 R\,\Fw((JH)^{\top}, \nabla \sin^2\theta)\leq 0
\end{equation}
then  $F$ is  either a complex or a Lagrangian submanifold. 
This is the case when $F$ has constant K\"{a}hler angle.
\end{Th}
\begin{Cr} Let $n=2$, $R<0$,  and $F:M\ra N$ be a closed submanifold
with parallel mean curvature and equal K\"{a}hler angles. 
If $\|H\|^2 \geq -\frac{R}{8}\sin^2\theta$, then $F$ is either a complex
or a Lagrangian submanifold.
\end{Cr}
\begin{Th} Assume $M$ is closed, $n\geq 3$,
and $F:M\ra N$ is an immersion with equal K\"{a}hler angles.\\[2mm]
$(A)$~~ If   $R< 0$, and if 
$ \delta\Fw((JH)^{\top})\geq 0 $,
then  $F$  is either complex or Lagrangian.\\[1mm]
$(B)$~~ If $R=0$,  and  if $\delta\Fw((JH)^{\top})\geq 0 $, then the
K\"{a}hler angle is constant.\\[1mm]
$(C)$~~If  $F$  has  constant K\"{a}hler angle and
 $R\neq 0$, then  $F$ is either complex or Lagrangian.
\end{Th}
In case $n=1 $ we obtain:
\begin{Pp} If $M$ is a closed surface and $N$ is a non Ricci-flat 
K\"{a}hler-Einstein surface, 
then any immersion $F:M\ra N$ either has  complex or  Lagrangian
points. In particular, if $F$ has constant  K\"{a}hler
angle, then  $F$ is either a complex or a Lagrangian submanifold.
\end{Pp}
This generalizes a result in [M-U], for
compact surfaces immersed with constant K\"{a}hler angle 
(and so orientable, if not Lagrangian) into $\Co\!I\!\!P^2$.\\[5mm]
For $M$ not  necessarily compact we have  the following proposition:
\begin{Pp} If $F:M\ra N$ is an immersion with constant equal K\"{a}hler
angle $\theta$ and with  parallel mean curvature, then:\\
(1)~ If $R=0$, $F$ is either   Lagrangian or minimal.\\[1mm]
(2)~ If $R>0$, $F$ is either   Lagrangian or complex.\\[1mm]
(3)~ If $R<0$, $F$ is either Lagrangian, or $\|H\|^2 = -\frac{\sin^2\theta}{4n}R.$\\
(4)~ If $H=0$, then $R=0$ or $F$ is either Lagrangian or complex.
\end{Pp}
Note that (4) of the above proposition is an improvement of
Theorem 1.3 of [S-V,2], for, compactness is not required now.
We also observe that from Corollary 1.1, if  $n=2$ and
$M$ were closed, that later case of (3) implies as well $F$ to be
complex or Lagrangian. Compactness of $M$ is a much more 
restrictive condition.
In [K-Z] it is shown that, if $n=1$ and $N$ is a complex space 
form of constant holomorphic sectional
curvature $4\rho$ and $M$ is a surface of non-zero parallel 
mean curvature and constant K\"{a}hler angle, 
 then either $F$ is Lagrangian and $M$ is flat,  or $\sin\theta=-\sqrt{
\frac{8}{9}}$, $\rho=-\frac{3}{4}\|H\|^{2}$ and $M$
has constant Gauss curvature $K=-\frac{\|H\|^2}{2}$. These values of 
$\theta$ and $\rho$ ($R=6\rho$) are according to our relation 
in (3) of Proposition 1.2. Chen in [Che,2] and [Che-T,2] shows explicitly
all possible examples of such (non-compact) surfaces of the 2-dimensional complex hyperbolic spaces.  In [K-Z]  it is also shown all  examples
of  surfaces immersed into $~\Co\! I\!\!H^2$  with non-zero 
parallel mean curvature  
and non-constant  K\"{a}hler angle.  In case (1), if $F$ is not minimal,
then $(JH)^{\top}$ defines a global nonzero parallel vector field on $M$
(see Proposition 3.6 of section 3).
\begin{Th} Let   $F$ be a closed surface immersed with
parallel mean curvature into a non Ricci flat K\"{a}hler-Einstein surface . 
If $F$ has  no complex points and if 
$\frac{\Fw}{Vol_{M}}\geq 0$ (or $\leq 0$) on all $M$, 
then  $F$ is Lagrangian. If $F$ has no Lagrangian points, then $F$ is
minimal.
\end{Th}

\section{Some formulas on the  K\"{a}hler angles}
\setcounter{Th}{0}
\setcounter{Pp}{0}
\setcounter{Cr} {0}
\setcounter{Lm} {0}
\setcounter{Def} {0}
\setcounter{equation} {0}
On $M$ we take the induced metric $g_{M}=F^{*}g$, that we also denote by $\langle,\rangle$.  We denote by  $\lnab{}$
both Levi-Civita connections of $M$ and  $N$, and by
$\lnab{X}dF(Y)=\lnab{}dF(X,Y)$
the second fundamental form of $F$, a symmetric tensor on $M$ with values on
the normal bundle $NM=(dF(TM))^{\bot}$ of $F$.
The mean curvature of $F$ is given by
$H=\frac{1}{2n}trace\lnab{}dF$.
At  each point $p\in M$, let 
$\{X_{\al}, Y_{\al}\}_{1\leq \al\leq n}$ be a $g_{M}$-orthonormal
basis of eigenvectors of $\Fw$. On that basis,  $\Fw$ is a $2n\times 2n$
block matrix
\[
F^{*}\omega=\bigoplus_{0\leq \al\leq n}  \left[ \begin{array}{cc}
             0 &  -\cos\theta_{\al} \\
             \cos\theta_{\al} &  0 \end{array} \right], 
\]
where $\cos\theta_{1}\geq \cos\theta_{2}\geq \ldots \geq \cos\theta_{n}
\geq 0$, are the corresponding eigenvalues
 ordered in decreasing way.  The angles $\{\theta_{\al}\}_{1\leq \al\leq n}$
are the K\"{a}hler angles of $F$ at $p$.
We identify the two form $\Fw$ with 
 the skew-symmetric operator of $T_{p}M$,
$(\Fw)^{\sharp}:T_{p}M\ra T_{p}M$,
using the musical isomorphism with respect to $g_{M}$, that is, $g_{M}((F^{*}\omega)^{\sharp}(X),Y)$ $=F^{*}\omega(X,Y)$, and  we take its  polar decomposition, 
$(F^{*}\omega)^{\sharp} =  |(F^{*}\omega)^{\sharp}| \, \Jw$,
where $J_{\omega}:T_{p}M\ra T_{p}M$ is a  partial isometry with the same 
kernel ${\cal K}_{\omega}$ as  of $F^{*}w$, and where
 $|(F^{*}\omega)^{\sharp}|=\sqrt{-(F^{*}\omega)^{\sharp 2} }$. 
On ${\cal K}^{\bot}_{\omega}$, the orthogonal complement of $\kw$ in 
$T_{p}M$, $J_{\omega}:{\cal K}^{\bot}_{\omega}\ra{\cal K}^{\bot}_{\omega}$ 
defines a $g_{M}$-orthogonal complex structure.     
On a open set without complex directions, that is $\cos\theta_{\al}<1$ $\forall\al$, we consider the locally Lipschitz map
\[\kappa=\sum_{1\leq \al\leq n}\log \left(\frac{1+\cos\theta_{\al}}{1-\cos\theta_{\al}}\right).\]
For each $0\leq k\leq n$, this map is  smooth on the largest  open set 
$\Omega^{0}_{2k}$, where 
$\Fw$ has constant rank $2k$.
  On a neighbourhood of a point $p_{0}\in
\Omega^{0}_{2k}$, we may take
 $\{X_{\al}, Y_{\al}\}_{1\leq \al\leq n}$  a smooth local 
$g_{M}$-orthonormal frame of $M$,  with $Y_{\al}=J_{\omega}X_{\al}$ for $\al\leq k$, and where $\{X_{\al},Y_{\al}\}_{\al\geq k+1}$  is 
any $g_{M}$-orthonormal frame of ${\cal K}_{\omega}$. 
Moreover, we may assume that this frame 
diagonalizes $F^{*}\omega$ at $p_{0}$. 
Following the computations of the appendix in [S-V,2],  
without requiring now minimality, we see that the components of the mean
curvature of  $F$  appear three times in the formula for $\triangle \kappa$.
Namely, when we compute $(5.9)$ and $(5.10)$ of [S-V,2], we get respectively, the extra terms $ig(\frac{n}{2}\lnab{\mu}H, JdF(\bmu))$ and 
$-ig(\frac{n}{2}\lnab{\bmu}H, JdF(\mu))$, and when  we sum 
$\sum_{\be}-R^M(\mu,\bbe,\be,\bmu)-R^M(\bmu,\bbe,\be,\mu)$ we obtain the
extra term $ng(H,\lnab{\mu}dF(\bmu))$. Then, we have to add in the final expression for
$\sum_{\be} Hess \tilde{g}_{\mu\bmu}(\be,\bbe)$ of Lemma 5.4 of [S-V,2]
the expression $\sum_{\be}ig(\frac{n}{2}\lnab{\mu}H, JdF(\bmu))
-ig(\frac{n}{2}\lnab{\bmu}H, JdF(\mu)) +\cos\theta_{\mu}ng(H,\lnab{\mu}dF(\bmu))$. Introducing these extra terms
in the  term $\sum_{\be,\mu}\frac{32}{\sin^2\theta_{\mu}}Hess
\tilde{g}_{\mu\bmu}(\be,\bbe)$ of $(5.7)$ of [S-V,2], we obtain
our  more general formula 
for $\triangle \kappa$:
\begin{Pp} For any immersion $F$, at a point $p_0$ on a open set where $\Fw$ has constant rank $2k$ and no complex directions, we have
\begin{eqnarray}
\triangle \kappa
&=&4i\sum_{\be} Ricci^{N}(JdF({\be}),dF({\bbe}))\\[-2mm]
&&+\sum_{\be,\mu}\!\frac{32}{\sin^{2}\theta_{\mu}}Im \La{(}\!
R^{N}(dF({\be}),dF({\mu}),dF({\bbe}),
JdF({\bmu})\!+\!i\cos\theta_{\mu}dF(\bmu))\!\La{)}\non\\[-2mm]
&&-\sum_{\be,\mu,\rho}\!\!\frac{64(\cos\theta_{\mu}
\!+\!\cos\theta_{\rho})}{\sin^{2}\theta_{\mu}\sin^{2}\theta_{\rho}}
Re\la{(}\gdf{\be}{\mu}{\brho}\gdf{\bbe}{\rho}{\bmu}\!\la{)}
\non \\[-2mm]
&& +\sum_{\be,\mu,\rho}\!\!\!
\frac{32(\cos\theta_{\rho}\!\!-\!\cos\theta_{\mu})}
{\sin^{2}\theta_{\mu}\sin^{2}\theta_{\rho}}\;(|\gdf{\be}{\mu}{\rho}|^{2}
\!\!+\!|\gdf{\bbe}{\mu}{\rho}|^{2})\non \\[-1mm]
&&+\sum_{\be,\mu,\rho}\frac{32(\cos\theta_{\mu}+\cos\theta_{\rho})}
{\sin^{2}\theta_{\mu}}\,\LA{(}
|\langle\lnab{{\be}}\mu,{\rho}\rangle|^{2} +
|\langle\lnab{{\bbe}}\mu,{\rho}\rangle|^{2}\LA{)}\non \\
&&+\!\!\sum_{\mu}\!\!\frac{8n}{\sin^{2}\theta_{\mu}}\LA{(}ig\la{(}
\lnab{\mu}H,\!JdF(\bmu)\la{)}
\!-\!ig\la{(}\lnab{\bmu}H,\!JdF(\mu)\la{)}\!+\!2\cos\theta_{\mu}
g( H,\lnab{\mu}dF(\bmu))\LA{)}\non
\end{eqnarray}
where $``\al"=Z_{\al}=\frac{X_{\al}-i Y_{\al}}{2}$ and $``\bar{\al}"=
\overline{Z_{\al}}$.
\end{Pp}

Projecting $JH$  on $dF(TM)$, we define a vector field $(JH)^{\top}$
on $M$,  and we denote by $((JH)^{\top})^{\flat}$ the corresponding 1-form, $((JH)^{\top})^{\flat}(X)$ $=g_{M}((JH)^{\top},X)$ $=g(JH,dF(X))$.
 If $F$ is a Lagrangian immersion, the above formula on $\triangle\kappa$
leads to a well-known result: 
\begin{Cr}{\bf ([W,2])} If $F$ is a Lagrangian immersion, then $((JH)^{\top})^{\flat}$ is
a closed 1-form on $M$.
\end{Cr}
A proof of this corollary  will be given in section 3.
The formula (2.1) is considerably simplified when
$F$ is an immersion with equal K\"{a}hler angles.
Now we recall the Weitzenb\"{o}ck formula for $\Fw$, that we used in 
[S-V,2]
\begin{equation}
\frac{1}{2}\triangle \|\Fw\|^{2} =-\langle \triangle \Fw,\Fw\rangle
+\|\lnab{}\Fw\|^{2}+\langle S\Fw, \Fw\rangle,
\end{equation}
where $\langle,\rangle$ denotes the Hilbert-Schmidt inner product for 2-forms,
and $S$ is the Ricci operator of $\bigwedge^{2}T^{*}M$,
and $\triangle = d\delta +\delta d$ is the the Laplacian operator on
forms.  $\Fw$ is a closed 2-form. 
 If it is also co-closed, that is $\delta\Fw=0$,
then it is harmonic. If $M$ is compact,
\begin{equation}
\int_{M}\langle\triangle\Fw, \Fw\rangle Vol_{M}
= \int_{M}\|\delta\Fw\|^{2}Vol_{M}
\end{equation}
We will use this formula when $F$ has equal K\"{a}hler angles.

\section{Immersions with equal  K\"{a}hler angles}
\setcounter{Th}{0}
\setcounter{Pp}{0}
\setcounter{Cr} {0}
\setcounter{Lm} {0}
\setcounter{Def} {0}
\setcounter{equation} {0}
In this section we recall some formulas for immersions with equal K\"{a}hler angles.  $F$ is said to have equal K\"{a}hler angles, if all the angles are
equal,  $\theta_{\al}=\theta$ $\forall \al$.  In this case, 
 $(\Fw)^{\sharp}=\cos\theta \Jw$,
 and $\Jw$ is a smooth almost complex structure away from
the set of Lagrangian points
${\cal L}=\{p\in M: \cos\theta(p)=0\}$. 
Let ${\cal L}^{0}$ denote the largest open set of ${\cal L}$, 
 ${\cal C}=\{p\in M: \cos\theta(p)=1\}$  the set of complex points, 
and ${\cal C}^{0}$ its largest open set.
Recall that $\cos^{2}\theta$ is smooth on all $M$, while
$\cos\theta$ is only locally Lipschitz on $M$, but smooth on
${\cal L}^{0}\cup (M\sim {\cal L})$. For
immersions with equal K\"{a}hler angles,  any local frame of the form
$\{X_{\al}, Y_{\al}=\Jw X_{\al}\}_{1\leq \al\leq n}$ diagonalizes
$\Fw$ on the whole set where it is defined. We use the letters $\al, \be,\mu,\ldots$
to range on the set $\{1,\ldots, n\}$ and the letters $j,k,\ldots$
to range on $\{1,\ldots, 2n\}$.
As in the previous section, we denote by $``\al"=Z_{\al}=\frac{X_{\al}-
i Y_{\al}}{2}$ and $``\bar{\al}"=\overline{Z_{\al}}=\frac{X_{\al}+
i Y_{\al}}{2}$, defining local frames on the complexifyied tangent space
of $M$.

 On tensors and forms we use the Hilbert-Shmidt inner product.
We denote by $\delta$ the divergence operator on (vector valued)
forms, and by
$div_{M}$ the divergence operator on vector fields over $M$.
The   $(1,1)$-part of
$\lnab{}dF$ with respect to  $\Jw$,  is given by
$(\lnab{}dF)^{(1,1)}(X,Y)=$ $\ha (\lnab{}dF(X,Y) + \lnab{}dF(\Jw X,\Jw Y))$.
This tensor is defined away from Lagrangian points, and it vanish
on ${\cal C}^{0}$, for, on that set, $F$ is a complex
submanifold of $N$, and $\Jw$ is the induced complex structure.
\begin{Pp}     {\bf [S-V,2]} On $(M\sim{\cal L})\cup {\cal L}^{0}$, 
\begin{eqnarray*}
\|\Fw\|^{2} & = & n\cos^{2}\theta\\
\|\lnab{}\Fw\|^{2} &=& n\|\nabla\cos\theta\|^{2}+
\frac{1}{2}\cos^{2}\theta\|\lnab{}\Jw\|^{2}\\
\delta(\Fw)^{\sharp} &=& (\delta\Fw)^{\sharp}
 =(n-2)\Jw(\nabla\cos\theta) \\
\|\delta\Fw\|^{2} &=& (n-2)^{2}\|\nabla\cos\theta\|^{2}\\
\cos\theta\delta\Jw &= &(n-1)\Jw(\nabla\cos\theta)
\end{eqnarray*}
and on $ \La{(}M\sim ({\cal L}\cup{\cal C})\La{)}\cup {\cal L}^{0}
\cup {\cal C}^{0}$,\\[-4mm]
\[(1\!-\!n)\nabla\sin^{2}\theta =
 16\cos\theta\, Re\LA{(} i\sum_{\be,\mu}
\La{(}\gdf{\bmu}{\mu}{\be}\!-\!
\gdf{\bmu}{\be}{\mu}\La{)}\bbe\LA{)}.
\]
In particular, for $n\neq 2$, $\Jw(\nabla\cos\theta)$, $\|\nabla\cos\theta\|^{2}$,
$\cos^{2}\theta\|\lnab{}\Jw\|^{2}$, and $\cos\theta\,\delta\Jw$
can be smoothly extended  to all $M$. Furthermore, for $n\geq 2$, there is a constant
$C>0$ such that on $M$, $\|\nabla\sin^2\theta\|^2 \leq C
\cos^2\theta€~\sin^2\theta~ \|(\lnab{}dF)^{(1,1)}\|^2$.
\end{Pp}
The  estimate on $\|\nabla\sin^2\theta\|^2$  given above follows from 
the expression on $(1-n)\nabla \sin^2\theta $ 
and the following  explanation. From Schwarz inequality,~
$|g(\lnab{X}dF(Y),JdF(Z))|= |g(\lnab{X}dF(Y),\Phi(Z))|\leq 
\|\lnab{X}dF(Y)\|\, \|\Phi(Z)\|$,
where $\Phi(Z)=(JdF(Z))^{\bot}$, and $(~)^{\bot}$ denotes the orthogonal
projection onto the normal bundle. 
 But (cf [S-V,2]) 
$JdF(Z) = \Phi(Z) + dF((\Fw)^{\sharp}(Z))$. 
An elementary computation
shows that 
\[
\|\Phi(Z)\|^2 = g\La{(}JdF(Z)\!-\!dF((\Fw)^{\sharp}(Z)),
JdF(Z)\!-\!dF((\Fw)^{\sharp}(Z))\La{)}
=\sin^2\theta\, \|Z\|^{2}
\]
Obviously the formula on $\nabla\sin^2\theta$ as well the estimate on $\|\nabla\sin^2\theta\|^2$,  are still valid on all
complex and Lagrangian points, since those points are critical points
for $ \sin^2\theta$, and at complex points $JdF(TM)\subset TM$ . Also 
\begin{Cr} If $n=2$, $\Fw$ is an harmonic 2-form. If $n\neq 2$,
$\Fw$ is co-closed iff $\theta$ is constant. For any $n\geq 2$, 
if $(M\sim {\cal L}, \Jw, g_{M})$ is K\"{a}hler, 
then $\theta$ is constant and $\Fw$ is parallel.
\end{Cr}
Following chapter $4$ of [S-V,2]
and  using the new expression for $\triangle\kappa$ of Proposition $2.1$,
with the extra  terms involving the mean curvature $H$, 
 and noting 
that  now  both $(4.4)$ and $(4.7)+(4.5)$  of [S-V,2] have extra terms  involving  $H$,  we obtain:
\begin{Pp} Away from complex and Lagrangian points, 
\begin{eqnarray*}
\lefteqn{
\triangle\kappa=}\\[-1mm]
&=&\!\!\!\!  \cos\theta \La{(}\! -2nR
+\frac{32}{\sin^{2}\theta}\sum_{\be,\mu}\!\! R^{M}(\be,\mu,\bbe,\bmu)
+\frac{1}{\sin^{2}\theta}\|\lnab{}\Jw\|^{2}
+\frac{8(n\!-\!1)}{\sin^{4}\theta}\|\nabla\cos\theta\|^{2}\La{)}\\[-1mm]
&&-\frac{16n}{\sin^{4}\theta}\cos\theta\sum_{\be}d\cos\theta\La{(}
i g (H, JdF(\be))\bbe-ig( H, JdF(\bbe))\be\La{)}\\[-1mm]
&&+\frac{8n}{\sin^{2}\theta}\sum_{\mu}\La{(}
ig( \lnab{\mu}H,JdF(\bmu))-ig( \lnab{\bmu}H,JdF(\mu))
\La{)}.
\end{eqnarray*}
\end{Pp}
Let us denote by $\lnabo{}$ the usual connection in the normal bundle,
and denote by $(JH)^{\top}$ the vector field of $M$ given by
\[g_{M} ((JH)^{\top}, X)= g(JH, dF(X))~~~~~\forall X\in TM.\\[3mm]\]
\begin{Lm} $\forall X,Y\in T_{p}M$,
\[\begin{array}{lrcl}
(i)&g(\lnab{X}H,JdF(Y)) &=& -\langle
\lnab{X}(JH)^{\top}, Y\rangle - g( H, J\lnab{X}dF(Y))~~~~~~~~~~~~~~(\mbox{on}~M)\\
& &=& -g(H,\lnab{X}dF ((\Fw)^{\sharp}(Y))) + g(\lnabo{X}H,JdF(Y))~~~~~(\mbox{on}~M)
\\[4mm]
(ii)&(\ha \Jw( (JH)^{\top})&=&\sum_{\be}i g(H, 
JdF(\be))\bbe-i g( H, JdF(\bbe))\be~~~~~~~~~~(\mbox{on~}M\sim{\cal L})
\end{array}\]
\begin{eqnarray*}
\!\!\!\!\!\!\!\!\!\!\!\!\!\!\!\!\!\!\!\!
(iii)~~~~~~~~~\lefteqn{\sum_{\mu}2i g( \lnab{\mu}H,\!JdF(\bmu))
-2ig( \lnab{\bmu}H,\!JdF(\mu))=}\\[-2mm]
&=&\sum_{\mu}  4 Im\langle \lnab{\mu} (JH)^{\top}, \bmu\rangle= 
-\sum_{\mu}2i d((JH)^{\top})^{\flat}(\mu, \bmu)~~~~~~~(\mbox{on}~M)\\[-2mm]
 &=& -2n\cos\theta\|H\|^2 - 4 \sum_{\mu} Im \La{(} g(\lnabo{\mu}H,JdF(\bmu))\La{)}~~~~~~~~~~~(\mbox{on}~M)\\[-1mm]
 &=&- div_{M}\La{(}\!\Jw((JH)^{\!\top})\La{)}
+ \langle (JH)^{\!\top}\!,\delta\Jw \rangle ~~~~~~~~~~~~~~~~~~~
(\mbox{on~}M\sim{\cal L}).\\[-3mm]
\end{eqnarray*}
(iv)~~~~~~~$div_{M}((JH)^{\top})  = \sum_{\mu}-4 Re\La{(} g( \lnabo{\mu}H, JdF(\bmu))\La{)}
~~~~~~~~~~~~~~~~~~(\mbox{on~}M)$.\\
\end{Lm}
\em Proof. \em
Assume that $\lnab{}~Y(p)=0$. Then we have at the point $p$
\begin{eqnarray*}
g( \lnab{X} H, JdF(Y)) &=& d \La{(}g( H, JdF(Y))\La{)} (X)
- g( H,  \lnab{X} (JdF(Y)))\\
&=& -d \langle (JH)^{\top}, Y\rangle (X)
- g( H, J \lnab{X}dF(Y) ) \\
&=& - \langle \lnab{X}(JH)^{\top}, Y \rangle
- g( H, J \lnab{X}dF(Y)).
\end{eqnarray*}
On the other hand,  from $JdF(Y)= dF( (\Fw)^{\sharp}(Y)) + (JdF(Y))^{\bot}$, we
get the second equality of $(i)$.
 For $p\in M\sim{\cal L}$, since $\Jw\be=i\be$, and  $\Jw\bbe=-i\bbe$, 
\begin{eqnarray*}
\lefteqn{\sum_{\be}i g( H, JdF(\be))  \bbe -i g( H, JdF(\bbe)) \be=}\\[-2mm]
&&\begin{array}{rcl}
&=& \sum_{\be}~g(H, JdF(\Jw\be)) \bbe + g( H, JdF(\Jw\bbe)) \be\\[2mm]
&=&\sum_{\be}~- g( JH, dF(\Jw\be)) \bbe - g( JH, dF(\Jw\bbe)) \be\\[2mm]
&=&\sum_{\be}~-\langle (JH)^{\top}, \Jw\be\rangle \bbe -\langle (JH)^{\top}, 
\Jw\bbe\rangle \be\\[2mm]
&=&\sum_{\be}~\langle \Jw((JH)^{\top}), \be\rangle \bbe +
\langle \Jw((JH)^{\top}), 
\bbe\rangle \be\\[2mm]
& =&\ha\Jw ((JH)^{\top}),
\end{array}
\end{eqnarray*}
and $(ii)$ is proved.
From the first equality of  $(i)$, 
\begin{eqnarray*}
\lefteqn{\sum_{\mu} i g(\lnab{\mu}H, JdF(\bmu)) 
-i g( \lnab{\bmu}H, JdF(\mu))=}\\[-1mm]
&=&  \sum_{\mu}-i\langle \lnab{\mu}(JH)^{\top}, \bmu\rangle 
+i\langle \lnab{\bmu}(JH)^{\top}, \mu\rangle=\sum_{\mu}2 Im \La{(}\langle \lnab{\mu}(JH)^{\top}, \bmu\rangle\La{)}  \\[-1mm]
&=&\sum_{\mu}-id((JH)^{\top})^{\flat}(\mu,\bmu).
\end{eqnarray*}
On the other hand, from second equality of $(i)$
\begin{eqnarray*}
\sum_{\mu}  g(\lnab{\mu}H, JdF(\bmu))&=& 
\sum_{\mu}-g(H,\lnab{\mu}dF(\cos\theta \Jw(\bmu)))
 +g(\lnabo{\mu}H,JdF(\bmu))\\
&=&\frac{ni}{2} \cos\theta\, g(H,H) +\sum_{\mu}g(\lnabo{\mu}H,JdF(\bmu)).
\end{eqnarray*}
Hence
\begin{eqnarray*}
\lefteqn{\sum_{\mu} i g(\lnab{\mu}H, JdF(\bmu)) 
-i g( \lnab{\bmu}H, JdF(\mu))=}\\[-1mm]
&=& -n\cos\theta\|H\|^2 -\sum_{\mu}2Im\La{(} 
g(\lnabo{\mu}H,JdF(\bmu))\La{)}.
\end{eqnarray*}
Similarly, from $div_{M}((JH)^{\top})=\sum_{\mu} 2\langle 
\lnab{\mu}(JH)^{\top}, \bmu\rangle + 2\langle 
\lnab{\bmu}(JH)^{\top}, \mu\rangle$ and $(i)$ we get $(iv)$.\\
Finally, 
 using the symmetry of $\lnab{}dF$ and that 
$\langle \lnab{Z}\Jw (X), Y\rangle
=-\langle \lnab{Z}\Jw (Y), X\rangle$ (cf. [S-V,2])
\begin{eqnarray*}
\lefteqn{\sum_{\mu} i g(\lnab{\mu}H, JdF(\bmu)) 
-i g( \lnab{\bmu}H, JdF(\mu))=}\\
&&\!\!\!\!\!\!\!\!\!\!\!\!\!\begin{array}{cl}
&= \sum_{\mu} \langle \lnab{\mu}(JH)^{\top}, \Jw(\bmu)\rangle 
+\langle \lnab{\bmu}(JH)^{\top}, \Jw(\mu)\rangle \\[2mm]
&= \sum_{\mu}
-\langle \Jw(\lnab{\mu}(JH)^{\top}), \bmu\rangle 
-\langle \Jw(\lnab{\bmu}(JH)^{\top}),\mu\rangle \\[2mm]
&= \sum_{\mu}\!\! -\langle\lnab{\mu} (\Jw(JH)^{\top}\!)
-\!\lnab{\mu}\Jw ((JH)^{\top}\!), \bmu\rangle 
-\langle \lnab{\bmu}(\Jw(JH)^{\top}\!)-\lnab{\bmu}\Jw((JH)^{\top}\!), 
\mu\rangle \\[2mm]
&=-\ha div_{M} (\Jw(JH)^{\top}) +\sum_{\mu}\langle \lnab{\mu}
\Jw ( (JH)^{\top}), \bmu\rangle + \langle \lnab{\bmu}
\Jw ( (JH)^{\top}), \mu\rangle\\[2mm]
&=-\ha div_{M} (\Jw(JH)^{\top}) +\sum_{\mu} -\langle (JH)^{\top},
\lnab{\mu}\Jw (  \bmu)\rangle - \langle  (JH)^{\top}, 
\lnab{\bmu}\Jw ( \mu)\rangle\\[2mm]
&=-\ha div (\Jw(JH)^{\top}) +\langle (JH)^{\top},
\ha\delta\Jw \rangle.~~~~~~~~~~~\qed\end{array}
\end{eqnarray*}\\[5mm]
Using $div(fX)=fdiv(X)+df(X)$, with $f=\frac{1}{\sin^{2}\theta}$,
and $X=\Jw( (JH)^{\top})$, and that $2\cos\theta d\cos\theta=
d\cos^{2}\theta=-d\sin^{2}\theta$,  we obtain applying Lemma $3.1$
to Proposition 3.2
\begin{Pp}
Away from complex and Lagrangian points
\begin{eqnarray*}
\triangle\kappa &=& {\sm \cos\theta \La{(} -2nR
+\frac{32}{\sin^{2}\theta}\sum_{\be,\mu} R^{M}(\be,\mu,\bbe,\bmu)
+\frac{1}{\sin^{2}\theta}\|\lnab{}\Jw\|^{2}
+\frac{8(n-1)}{\sin^{4}\theta}\|\nabla\cos\theta\|^{2}~\La{)}}\\
&&-div_{M}\left(\Jw\LA{(}\frac{4n(JH)^{\top}}{\sin^{2}\theta}\LA{)}
\right) + g_{M}\LA{(}\delta\Jw, \frac{4n(JH)^{\top}}{\sin^{2}\theta}\LA{)}.
\end{eqnarray*}
\end{Pp}
If $n=1$ then $(M,\Jw,g)$ is a K\"{a}hler manifold (away from
 Lagrangian points),
and so, $\delta \Jw=\lnab{}\Jw=0$. Obviously the curvature term on $M$ 
in the expression of $\triangle\kappa$ vanishes. Then, $\triangle\kappa$ reduces to:
\begin{Cr}
If $n=1$, away from complex and Lagrangian points
\begin{equation}
\triangle\kappa = -2R\cos\theta -4div_{M}\left(
\Jw\LA{(}\frac{(JH)^{\top}}{\sin^{2}\theta}\LA{)}.
\right) 
\end{equation}
\end{Cr}
Now we compute $\triangle\cos^2\theta$ from $\triangle\kappa$ of
Proposition 3.3 and applying  Proposition 3.1,
following step by step the proof of Proposition 4.2 of [S-V, 2].
Recall that,  if $F$ has equal K\"{a}hler angles  at $p$,
then, at $p$ (cf.[S-V,2])
 \[\langle S\Fw,\Fw\rangle =16\cos^{2}\theta\sum_{\rho,\mu}
R^{M}(\rho,\mu,\brho,\bmu),\]
where $S\Fw$ is the Ricci operator applied to $\Fw$,  appearing in the
 Weitzenb\"{o}ck formula $(2.2)$. 
If $(M,\Jw,g_{M})$ is
K\"{a}hler in a neighbourhood of $p$, then 
$\langle S\Fw,\Fw\rangle=0$ at $p$.
\begin{Pp}
Away from complex and Lagrangian points:
\begin{eqnarray}
n\triangle\cos^{2}\theta &=&  -2n\sin^{2}\theta\cos^{2}\theta R + 
2\langle S\Fw, \Fw\rangle + 2\|\lnab{}\Fw\|^{2}\non \\
&& +4(n-2)\|\nabla \,|\sin\theta|~\|^{2}-4n \,div_{g_{M}}\left( (\Fw)^{\sharp}((JH)^{\top})\right) \non \\
&&-\frac{4n(2 +(n-4)\sin^2\theta)}{\sin^{2}\theta}\langle 
\nabla\cos\theta, \Jw( (JH)^{\top})\rangle
\end{eqnarray}
The last term $(3.2)$  can be written, for  $n=2$ as
\begin{equation}
(3.2)= 8\,\Fw ((JH)^{\top},\nabla \log\sin^2\theta)
\end{equation}
and for $n\geq 3$,
\begin{equation}
(3.2)= \frac{4n(2 +(n-4)\sin^2\theta)}{\sin^{2}\theta(n-2)}
\delta\Fw((JH)^{\top})
\end{equation}
\end{Pp}
The expressions in  (3.3) and (3.4) come from Proposition $3.1$ and 
the fact that
$(\Fw)^{\sharp}=\cos\theta\Jw$.\\[8mm]
\em  Remark 1.~\em  Let $\wt=\omega_{| NM}$ be the restriction of 
the K\"{a}hler form
$\omega$ to the  normal vector bundle $NM$, and
$\wt= |\wt|\Jt$ be its polar decomposition, when we identify it with a
skew-symmetric operator on the normal bundle, using the musical isomorphism.
Let $\cos\sigma_{1}\geq \cos\sigma_{2}\geq \ldots\geq \cos\sigma_{n}\geq 0 $
be the eigenvalues of $\wt$. 
The $\sigma_{\al}$ are the K\"{a}hler angles of $NM$. 
If $\{U_{\al},V_{\al}\}$ is an orthonormal basis of eigenvectors of $\wt$
at $p$, then 
$\wt = \sum_{\be}\cos\sigma_{\be} U_{*}^{\be}\wedge V_{*}^{\be}$.
  For each $p$,   $CD(F)=  \bigoplus _{\al:~\cos\theta_{\al}=1}
 span  \{ X_{\al}, Y_{\al}\}$  defines 
the vector subspace of complex directions, or equivalently, 
the largest $J$-complex vector subspace contained in $T_{p}M$. Similarly 
we define $CD(NM)$, the largest $J$-complex subspace of $NM$ at $p$.
Then  
\[\begin{array}{ccl}
\Fw&=& \omega_{|CD(F)} +\sum_{\cos\theta_{\al} <1
}\cos\theta_{\al} X_{*}^{\al}\wedge Y_{*}^{\al}\\[2mm]
\wt& =& \omega_{|CD(NM)} +
\sum_{\cos\sigma_{\al} <1}\cos\sigma_{\al} U_{*}^{\al}\wedge V_{*}^{\al}
\end{array}\]
We define the following morphisms between vector bundles of 
the same dimension $2n$, where $(~)^{\top}$  and $(~)^{\bot}$ 
denote the orthogonal
projection onto $TM$ and $NM$ respectively,
\[\begin{array}{cccccc}
 \Phi : TM & \ra  & NM ~~~~~~~~~~~~~~~~~ \Xi : NM & \ra & TM\\
         X &  \ra  & (JdF(X))^{\bot}~~~~~~~~~~~~~~ U & \ra & (JU)^{\top}
\end{array}\]
Then $\Phi^{-1}(0)=CD(F)$, $\Xi^{-1}(0)=CD(NM)$.
Note that $\forall X,Y\in TM$ and $\forall U,V\in NM$
\[\begin{array}{cl}
(JdF(X))^{\top}= dF((\Fw)^{\sharp} (X))&~~~(JU)^{\bot}= \wt (U),\\
\Phi(X)= JdF(X) - dF((\Fw)^{\sharp}(X))&~~~\Xi(U)= JU - \wt(U).\end{array}\]
A simple computation shows that, if $\cos\theta_{\al}\neq 1$, we may take
  ~$U_{\al}=\Phi(\frac{Y_{\al}}{\sin\theta_{\al}})$, and
$V_{\al}=\Phi(\frac{X_{\al}}{\sin\theta_{\al}})$. Moreover, $CD(NM)=CD(F)^{\bot}\cap NM$ and $dim~CD(F)= dim~CD(NM)$.
Then  $\wt$ and $\Fw$ have the same eigenvalues, that is $NM$ and
$F$ have the same K\"{a}hler angles.
We  also define $LD(F)=Ker~\Fw={\cal K}_{\omega}$, $LD(NM)=Ker~\wt$ the 
vector subspaces of Lagrangian directions of $F$ and $NM$ respectively.
Then  we have $J(LD(F))=LD(NM)$. Furthermore,
 $\Jt\circ \Phi= -\Phi\circ \Jw$,~ 
  $\Jw\circ \Xi = -\Xi\circ \Jt$,~ 
$ -\Xi\circ \Phi = Id_{TM} + ((\Fw)^{\sharp})^{2}$,~ 
$ -\Phi\circ \Xi =Id_{NM} + (\wt)^{2}$.
 Considering the Hilber-Smidt norms,
$\|\Phi\|^{2}=\|\Xi\|^{2}=2\sum_{\al}\sin^{2}\theta_{\al}$.
If $F$ has equal K\"{a}hler angles, $-\Xi\circ \Phi=\sin^{2}\theta Id_{TM}$,
 ~~$-\Phi\circ \Xi=\sin^{2}\theta Id_{NM}$, and
\[
g(\Phi(X),\Phi(Y))=\sin^2\theta \langle X,Y\rangle~~~~~~
\langle\Xi(U),\Xi(V)\rangle =\sin^2\theta\, g(U,V).
\]
If $F$ has equal K\"{a}hler angles, since $NM$ and $F$
 have the same K\"{a}hler angles, we see that,  at a point $p\in M$ 
such that $H\neq 0$,   $(JH)^{\top}=0$
iff $ p$ is a complex point of $F$. We also note that, from lemma 3.1(iv),
if $F$ has parallel mean curvature $(JH)^{\top}$ is divergence-free, or
equivalentely, $((JH)^{\top})^{\flat}$ is co-closed.\\[2mm]

In [S-V,2] we have defined non-negative isotropic scalar curvature, as
a less restrictive condition than  non-negative isotropic
sectional curvature of [Mi-Mo]. If such curvature condition on $M$
holds, then $\sum_{\rho,\mu}R^{M}(\rho,\mu,\brho,\bmu)\geq 0$,
where $\{\rho, \brho\}_{1\leq \rho\leq n}$ is the complex basis  
of $T_{p}^c M$ defined by a basis of eigenvectors of $\Fw$.
Hence, if $F$ has equal K\"{a}hler angles  
$\langle S\Fw,\Fw\rangle\geq 0$.
A simple  application of the  Weitzenb\"{o}k formula (2.2)  shows 
in next proposition,  that such curvature condition on $M$, implies
the angle must be constant. No minimality is required.

\begin{Pp}  ({\bf [S-V,2]}) Let $F$ be a  non-Lagrangian 
immersion with equal K\"{a}hler angles of a compact
 orientable $M$ with non-negative isotropic scalar curvature into a 
K\"{a}hler manifold $N$.
If $n=2$, $3$ or $4$, then $\theta$ is constant and $(M,\Jw,g_{M})$ is a K\"{a}hler manifold.   For any $n\geq 1$ and $\theta$ constant, 
$\Fw$ is parallel, that is,  $(M,\Jw,g_{M})$ is a K\"{a}hler manifold.
\end{Pp}
Finally, before we prove  Corollary 2.1, we state a more general proposition.
Let $F:M\ra N$ be an immersion with equal K\"{a}hler angles, and let 
$M'=\{p\in M:  H=0\}$ be the set of minimal points of $F$. On $M\sim {\cal C}$
 a 1-form is defined
\[ \sigma=\frac{2n}{\sin^2\theta} ((JH)^{\top})^{\flat} +
\frac{\delta \Fw}{\sin^2\theta}\]
Following the proof of [G], but now neither requiring $n=2$ nor
$\delta \Fw=0$, we obtain
\[\begin{array}{lcl}
\sigma(X) &=& -trace~ \frac{1}{\sin^2\theta}
 g(\lnab{}dF(\cdot, X), JdF(\cdot))\\[2mm]
d\sigma(X,Y) &=&
Ricci^{N}(JdF(X),dF(Y))=R\Fw(X,Y)\end{array}\]
We note that this form $\sigma$ is well known (see e.g [Br], [Che-M], [W,2]).
Now we have:\\[-2mm]
\begin{Pp} If $n=2$,
or if $n\geq 2$ and $\theta$ is constant,  then $\sigma=\frac{2n}{\sin^2\theta}((JH)^{\top})^{\flat}$
and  does not vanish on $M\sim (M'\cup {\cal C})$. Moreover,
 if $R=0$, then $d\sigma=0$. Thus, 
 if $\theta$ is constant $\neq 0$,  $\sigma\in H^{1}(M,\R{})$, and in
 particular, if $F$ has non-zero parallel mean curvature,
and $R=0$,  then $F$ is Lagrangian and  $\sigma$ is a non-zero 
parallel 1-form on $M$.\\[2mm]
For any immersion  with constant equal K\"{a}hler angles, the following
equalities hold
\[R\cos\theta\sin^2\theta = \sum_{\be}2d ((JH)^{\top})^{\flat}(X_{\be},Y_{\be})=-4n\cos\theta\|H\|^2 -\sum_{\mu}8 Im\La{(} g(\lnabo{\mu}H, JdF(\bmu))\La{)},\]\\[-5mm]
where $\{X_{\al},Y_{\al}\}$ is any basis of eigenvalues of $\Fw$.\\
\end{Pp}
\em Proof of Proposition 3.6 and Corollary 2.1. \em
We start by proving Corollary 2.1.
For a Lagrangian immersion, the formula on $\triangle\kappa$ (valid on
$\Omega^{0}_{0}$), 
reduces to
\[
 0=\triangle\kappa = \sum_{\mu,\be} 32 Im \La{(}
R^N (dF(\be), dF(\mu),dF(\bbe), JdF(\bmu))\La{)} -
\sum_{\mu}16n Im \La{(} g\la{(}
\lnab{\mu}H,JdF(\bmu)\la{)}\La{)}.
\]
Applying Codazzi equation to the curvature term
and noting that $JdF(TM)$ is the orthogonal complement of $dF(TM)$,
and that  $\sum_{\be}\lnab{\mu}\lnab{}dF (\be,\bbe) =\frac{n}{2}\lnabo{\mu}H$, we get\\[-3mm]
\begin{equation}
 0= \sum_{\be,\mu} Im ~\LA{(} g\La{(} \lnab{\be}\lnab{}dF (\mu,\bbe),
JdF(\bmu)\La{)}\LA{)}.
\end{equation}
Note that,  since $F$ is Lagrangian, we can choose arbitrarily the orthonormal 
frame $X_{\al},Y_{\al}$. Then we may assume they have zero covariant 
derivative at a given point $p$.
Since $F$ is a Lagrangian immersion 
$g(\lnab{}dF(\be,\bmu),JdF(\mu))=g(\lnab{}dF(\bmu,\mu),JdF(\be))$ (see e.g [S-V,2]).
Taking the derivative of  this equality at the point $p$ in the direction $\bbe$ we obtain
\begin{eqnarray*}
\lefteqn{g\La{(} \lnab{\bbe}\lnab{}dF (\be,\bmu),JdF(\mu)\La{)}
+g\La{(} \lnab{}dF (\be,\bmu),J\lnab{}dF(\bbe, \mu)\La{)}=}\\
&&g\La{(} \lnab{\bbe}\lnab{}dF (\bmu,\mu),JdF(\be)\La{)}
+g\La{(} \lnab{}dF (\bmu,\mu),J\lnab{}dF(\bbe, \be)\La{)}.
\end{eqnarray*}
Taking the summation on $\mu,\be$ and the imaginary part, we obtain
 from (3.5)
\[ \sum_{\be} Im~\La{(} g(\lnab{\bbe}H,JdF(\be) )\La{)}=\sum_{\be} Im~\La{(} g(\lnabo{\bbe}H,JdF(\be) )\La{)} =0.\]
From Lemma 3.1 we conclude,
\[\ha i\sum_{\be}d((JH)^{\top})^{\flat}(X_{\be},Y_{\be})=
-\sum_{\be}d((JH)^{\top})^{\flat}(\bbe,\be)=
\sum_{\be}-2i Im g_{M} (\lnab{\bbe}(JH)^{\top},\be)=
0.\]
From the arbitrarity of the orthonormal frame, we may interchange  $X_1$ by $-X_{1}$,  obtaining  $d((JH)^{\top})^{\flat}(X_{1},Y_{1})=0$. Hence $d((JH)^{\top})^{\flat}=0$.\\[5mm]
Now we prove Proposition 3.6.
The first part is an immediate conclusion from the expressions for $\sigma$,
$d\sigma$, and the fact that, under the above assumptions,  
$\delta \Fw=0$ (see Corollary 3.1), besides
 the  considerations on the zeroes  of $(JH)^{\top}$
in the previous remark. The conclusion that $F$ is Lagrangian
and $\sigma$ is parallel, under the
assumption of non-zero parallel mean curvature and $R=0$, comes from
the equalities stated in the proposition, which   we prove now,
and from Lemma  4.1 of  next section .
 It is obviously true if $\cos\theta= 1$,
that is for complex immersions, and it is true for $\cos\theta=0$, as we have seen above. Now, if $\cos\theta$ is constant and different
from $0$ or $1$, from Proposition 3.3,
\begin{eqnarray*}
0=\triangle\kappa &=&  \cos\theta \La{(} -2nR
+\frac{32}{\sin^{2}\theta}\sum_{\be,\mu} R^{M}(\be,\mu,\bbe,\bmu)
+\frac{1}{\sin^{2}\theta}\|\lnab{}\Jw\|^{2}
~\La{)}\\
&&-\frac{4n}{\sin^{2}\theta}div_{M}\left(\Jw\LA{(}(JH)^{\top}\LA{)}
\right) +\frac{4n}{\sin^{2}\theta} g\LA{(}\delta\Jw, (JH)^{\top}\LA{)}.
\end{eqnarray*}
Since $\Fw$ is harmonic   (see Corollary 3.1), 
 Weitzenb\"{o}ck formula (2.2) with $\theta$ constant reduces to
\[ 16\cos^2\theta \sum_{\be,\mu} R^{M}(\be,\mu,\bbe,\bmu)
=\langle S\Fw,\Fw\rangle=-\|\lnab{}\Fw\|^{2}=-\ha\cos^2\theta
\|\lnab{}\Jw\|^2\]
Thus,  from lemma 3.1
\begin{eqnarray*}
\ha R\cos\theta\sin^2\theta&=&-div_{M}\left(\Jw\LA{(}(JH)^{\top}\LA{)}
\right) +g_{M}\La{(}\delta\Jw, (JH)^{\top}\LA{)}\\[-1mm]
&=& -2n\cos\theta \|H\|^2 -4\sum_{\mu} Im\La{(} g(\lnabo{\mu} H,JdF(\bmu))\La{)}.~~~~\qed
\end{eqnarray*}

\section{Proofs of the main results}
\setcounter{Th}{0}
\setcounter{Pp}{0}
\setcounter{Cr} {0}
\setcounter{Lm} {0}
\setcounter{Def} {0}
\setcounter{equation} {0}
\em Proof of Proposition 1.1. \em Assume ${\cal C}\cup {\cal L}=\emptyset$.
Then the formula in Corollary 3.2 is valid on all $M$ with all 
maps involved smooth everywhere.  By applying
Stokes we get $\int_{M} R \cos\theta\, Vol_{M}=0$, where $\cos\theta>0$,
which is impossible if $R\neq 0$.\qed\\[6mm]
\em Proof of Proposition 1.2. \em Follows immediately from Proposition 3.6.
\qed\\[6mm]
\em Proof of Theorem 1.4. \em
In case $n=1$, $\Fw$ is a multiple of the volume element of $M$, that is
$\Fw=\cos\tilde{\theta}Vol_{M}$.
This $\tilde{\theta}$ is  the genuine definition of K\"{a}hler angle
given by Chern and Wolfson [Ch-W]. Our is just
$\cos\theta= |\cos\tilde{\theta}|$.
While $\cos\tilde{\theta}$ is smooth on all $M$, $\cos\theta$
may not be $C^{1}$ at Lagrangian points.
 But we see that the
formula $(3.1)$ is also valid on $M\sim {\cal L}\cup {\cal C}$
replacing $\cos\theta$ by $\cos\tilde{\theta}$ and the corresponding
replacement of $\kappa$ by $\tilde{\kappa}$, and  $\sin^2\theta$
by $\sin^2\tilde{\theta}$ and
$\Jw$ by $J_{M}$, the  natural 
$g_{M}$-orthogonal complex structure on $M$,  defining a K\"{a}hler 
structure. We denote this new formula by $(3.1)'$.
 Note that on  $M\sim {\cal L}$, $\Jw=\pm J_{M}$, 
the sign being $+$  or $-$ according to the sign of $\cos\tilde{\theta}$.
Hence a change of the sign of $\cos\tilde{\theta}$ will give a change
of sign on $\tilde{\kappa}$ and on $\Jw$ (w.r.t. $J_{M}$).
The formula $(3.1)'$ is  in fact also valid on  ${\cal L}^{0}$. 
To see this we use  the following lemma, as an immediate
consequence of Lemma $3.1~(i)$:
\begin{Lm} If $F:M^{2n}\ra N^{2n}$  is a  submanifold 
with parallel mean curvature, then $(JH)^{\top}$ is a parallel
vector field along ${\cal L}$, that is $\lnab{}(JH)^{\top}(p)=0$~
$\forall p\in {\cal L}$. 
\end{Lm}
Now it follows that $div_{M}( J_{M}((JH)^{\top}) )=0$ on ${\cal L}$.
Hence, the formula $(3.1)'$ on $\triangle \tilde{\kappa}$  
 is valid on ${\cal L}^0$, that is, at interior 
Lagrangian points. If we assume ${\cal C}=\emptyset$,  then $(3.1)'$
 is valid over all $M$, because now $\tilde{\kappa}$, 
$\cos\tilde{\theta}$, $J_{M}$,  and 
$\sin^2\tilde{\theta}$
are  smooth everywhere and ${\cal L}
\sim{\cal L}^{0}$ is a set of Lagrangian points with no interior.
Integrating and using Stokes, $2R\int_{M}\cos\tilde{\theta}
=0$. Hence if $\cos\tilde{\theta}$ is  non-negative 
or non-positive everywhere,  and if $R\neq 0$, then  $F$ is Lagrangian.
If $F$ has no Lagrangian points, from Lemma 3.1 $(iii)$, since 
$\delta\Jw=0$, 
\[div_{M}\La{(} \Jw(JH)^{\top}\La{)}=2\cos\theta\|H\|^2\]
is valid on $M$. Integration leads to $H=0$.
\qed.\\[5mm]
\em Proof of Theorem 1.2. \em
If $n=2$, using (3.3) in
the expression of $\triangle\cos^{2}\theta$ in Proposition 3.4, 
we get an expression that is smooth away from
complex points, and valid at interior Lagrangian points, and hence
on all $M\sim {\cal C}$. Then, following the same steps in the proofs of
 [S-V,2] chapter 4,  combining the formulae for $\triangle\cos^2\theta$ 
of Proposition 3.4 and the Weitzenb\"{o}k formula (2.2),  and applying
 Proposition 3.1, we get,
 away from complex  points
\begin{equation}
 sin^{2}\theta\cos^{2}\theta R =-2div_{M}((\Fw)^{\sharp}((JH)^{\top}))
+ 2\Fw((JH)^{\top}, \nabla\log \sin^2\theta)
\end{equation}
Set $P= sin^{2}\theta\cos^{2}\theta R +2 div_{M}((\Fw)^{\sharp}((JH)^{\top}))$.  This map is defined and smooth on all $M$ and vanishes on  ${\cal C}^{0}$. 
If $R>0$ (resp. $R<0$), and under the assumption $(1.1)$, we have from (4.1)
 that $P\leq 0$ (resp. $\geq 0$) on $M\sim {\cal C}$. Since the remaining set ${\cal C}\sim {\cal C}^0$
is a set of empty interior, then $P\leq 0$ (resp. $\geq 0$) is
valid on all $M$. In fact,  from Proposition 3.1, 
  $ |\Fw((JH)^{\top}, \nabla \sin^2\theta)|
\leq \sqrt{C}\cos^2 \theta\sin^2\theta\| H\|~\|(\lnab{}dF)^{(1,1)}\|$.
Since  $ (\lnab{}dF)^{(1,1)}$ vanishes on ${\cal C}^{0}$, and so also on $
\overline{{\cal C}^{0}}$, we can smoothly extend to zero
$\Fw((JH)^{\top}, \nabla \log\sin^2\theta)$ on $\overline{{\cal C}^0}$.
This we can also get from (4.1). Moreover, such equation tells us
we can smoothly extend the last term to all complex points, giving exactly
the value $2div_{M}((\Fw)^{\sharp}((JH)^{\top}))$ at those points.
 Integration of $P\leq 0$ (respectively $\geq 0$) and applying
Stokes, we  have
\[\int_{M}\sin^2\theta\cos^2\theta R Vol_{M}\leq 0 \mbox{~~~~(resp.~~}\geq 0)
\]
and conclude that   $F$ is either complex  or Lagrangian.~~~\qed \\[5mm]
\em Proof of Corollary 1.1. \em   Instead of using Stokes on
the term $div_{M}\La{(}(\Fw)^{\sharp} ( (JH)^{\top}) )\La{)}$, to make it disapear as we did in the proof of theorem 1.2, 
we develop it  into 
\begin{eqnarray*}
 div_{M}\La{(}(\Fw)^{\sharp} ( (JH)^{\top}) )\La{)}&=&
div_{M}\La{(}\cos\theta\Jw ( (JH)^{\top}) )\La{)}\\
&=& \cos\theta div_{M}\La{(}\Jw ( (JH)^{\top}) )\La{)} +
d\cos\theta \La{(}\Jw ( (JH)^{\top}) )\La{)},
\end{eqnarray*}
and use Lemma 3.1 to give, away from complex and Lagrangian points,
\begin{eqnarray*}
sin^{2}\theta\cos^{2}\theta R&=&-2\cos\theta div_{M}(\Jw((JH)^{\top}))
-2\langle \Jw( (JH)^{\top}), \nabla \cos\theta \rangle\\
&&+2\Fw((JH)^{\top}, \nabla\log \sin^2\theta)\\
&=& -8\cos^2\theta\| H\|^2+2\Fw((JH)^{\top}, \nabla\log \sin^2\theta).
\end{eqnarray*}
Hence, away from complex and Lagrangian points
\[ \sin^4\theta \cos^2\theta R + 8 \sin^2\theta \cos^2\theta\| H\|^2
=2\Fw((JH)^{\top}, \nabla \sin^2\theta).\]
Obviously, this equality also holds at Lagrangian and complex points,
for, those points are critical points for $\sin^2\theta$.
The corollary now follows immediately from Theorem 1.2. ~~~\qed\\[5mm]
\em Proof of Theorem 1.3.~ \em 
If $n\geq 3$ we  set 
\[P= n\triangle\cos^2\theta +
4n\, div_{M}((\Fw)^{\sharp}((JH)^{\top})) +2n\sin^2\theta\cos^2\theta R -2\|\lnab{}\Fw\|^2 -2\langle S\Fw, \Fw\rangle.\]
This map is defined on all $M$ and is smooth. From Proposition (3.4)
and using (3.4), on $M\sim {\cal C}$
\[P= \frac{4n(2+(n-4)\sin^2\theta)}{(n-2)\sin^{2}\theta}
\delta\Fw((JH)^{\top}) + 4(n-2)\|\nabla |\sin\theta|~\|^2\]
In $(A)$ and $(B)$, by assumption, $P\geq 0$ on $M\sim{\cal C}$,
because for $n\geq 3$, $ (2 + (n-4)\sin^2\theta) \geq 0$.
But on ${\cal C}^0$, $P=0$, for $(M,\Jw,g_{M})$ is
a complex submanifold, and so,
   $(JH)^{\top}=0$ and $\langle S\Fw, \Fw\rangle =0$.
Thus, $P\geq 0$ on all $M$. Integrating $P\geq 0$ on  $M$
we obtain using Stokes,  Weitzenb\"{o}ck formula (2.2), and (2.3)
\[ \int_{M}2n R\sin^2\theta\cos^2\theta Vol_{M}\geq \int_{M}2\|
\delta\Fw\|^2Vol_{M}.\]
Thus, if $R<0$ we conclude $F$ is either complex or Lagrangian,
and if $R=0$ we conclude that $\delta\Fw=0$, which implies, by 
Corollary 3.1, that $\theta$ is constant. This last reasoning
proves $(C)$ as well.\qed \\[7mm]
\em Remark 2. \em   In Theorem 1.3 we can replace the condition $\delta
\Fw ((JH)^{\top})\geq 0$ by a weaker condition 
\[\delta\Fw ((JH)^{\top})\geq -\frac{(n-2)^{2}}{4n(2+(n-4)\sin^2\theta)}
\|\nabla \cos^2\theta\|^2\]
 to achieve the same conclusion. This condition is sufficient to obtain $P\geq 0$ in the above proof. Then we can obtain for $n\geq 3$ 
a  corollary similar to Corollary 1.1,  by requiring 
\[ 4n^2 \cos^2\theta\|H\|^2 + n \sin^2\theta\cos^2\theta R -(n-2)^2
\|\nabla \cos\theta\|^2\geq -2n\delta\Fw ((JH)^{\top}).\]
{\Large \bf{References}}\\[2mm]
[Br]~~{\small R.L.\ Bryant,  \em Minimal Lagrangian submanifolds of 
K\"{a}hler-Einstein manifolds \em,   in \em Differential Geometry and Differential
Equations \em (Shangai, 1985), Lecture Notes is Math. {\bf 1255}, Springer,
Berlin, (1987), 1-12.}\\[1mm]
[Che,1]~~{\small B-Y.\ Chen, \em  Geometry of slant submanifolds, \em
Katholieke Universiteit Leuven, 1990}\\[1mm]
[Che,2]~~{\small B-Y.\ Chen,\em Special slant surfaces and a basic inequality,
\em Results Math.{\bf 33}(1998),65-78.}\\[1mm]
[Che-M]~~{\small B-Y.\ Chen \& J-M.\ Morvan, \em Cohomologie des 
sous-vari\'{e}t\'{e}s $\alpha$-oblique, \em C.R.\ Acad.\ Sci.\. Paris {\bf 314}
(1992),  931-934.}\\[1mm]
[Che-T,1]~~{\small B-Y.\ Chen  \& Y.\ Tazawa, \em Slant submanifolds in complex Euclidean  spaces, \em   Tokyo J. Math. {\bf 14} (1991) 101-120.}\\[1mm]
[Che-T,2]~~{\small B-Y.\ Chen  \& Y.\ Tazawa, \em Slant submanifolds of complex projective and complex hyperbolic spaces,\em  Glasgow Math.\ J. {\bf 42} (2000) 439-454.}\\[1mm]
[Ch-W]~~{\small S.S.\ Chern \& J.G.\ Wolfson, \em Minimal surfaces by moving frames, \em Amer.\  J.\  Math. {\bf 105} (1983), 59-83.}\\[1mm]
[D-S]~~{\small D.M.\ Duc \& I.M.C. Salavessa, \em Graphs with equal K\"{a}hler angles, \em in preparation.}\\[1mm]
[G]~~{\small A.\ Ghigi, \em  A generalization of Cayley submanifolds, \em
IMRN {\bf 15} (2000), 787-800.}\\[1mm]
[H-L]~~{\small F.R.\ Harvey \& H.B.\ Lawson Jr, \em Calibrated geometries,
\em Acta Math. {\bf 148} (1982), 47-157.}\\[1mm]
[K-Z]~~{\small K.\ Kenmotsu \& D.\ Zhou, 
\em The classification of the surfaces with parallel mean curvature
vector in two-dimensional complex space forms, \em  Amer.\ J.\ Math.
{\bf 122} (2000), 295--317.}\\[1mm]
[Mi-Mo]~~{\small M.J.\ Micallef \& J.D.\ Moore, \em
Minimal two-spheres and the topology of manifolds with positive
curvature on totally isotropic two-planes, \em Annals of Math. {\bf 127}
(1988), 199-227.}\\[1mm] 
[M-U] ~~{\small S.\ Montiel \& F.\ Urbano,
\em A Willmore functional for compact surfaces of 
complex projective plane, \em  preprint arXiv:math.DG/0002155}\\[1mm]
[S-V,1]~~{\small I.M.C.\ Salavessa \& G.\ Valli, \em
Broadly-Pluriminimal Submanifolds of K\"{a}hler-Einstein Manifolds, \em
 Yokohama Math.\ J. {\bf 48} (2001), 181-191.} \\[1mm]
[S-V,2]~~{\small I.M.C.\ Salavessa \& G.\ Valli, \em Minimal
submanifolds of K\"{a}hler-Einstein manifolds with equal K\"{a}hler angles,
\em To appear in Pacific J.\ Math. (a previous version: e-print no.\ math.DG/0002050).}\\[1mm]
[T]~~{\small  G.\ Tian, \em Gauge theory and calibrated geometry, I, \em
Annals of Math.{\bf 151} (2000), 193-268. }\\[1mm]
[W,1]~~{\small J.G.\ Wolfson, \em Minimal
Surfaces in K\"{a}hler Surfaces and Ricci Curvature, \em J.\ Diff.\ Geo.
{\bf 29} (1989), 281--294.}\\[1mm]
[W,2]~~{\small J.G.\ Wolfson, \em Minimal Lagrangian Diffeomorphisms
and the Monge-Amp\`{e}re Equation, \em J.\ Diff.\ Geom.
{\bf 46} (1997), 335-373.}\\[4mm]
 {\footnotesize Centro de F\'{\i}sica das
Interac\c{c}\~{o}es Fundamentais,\\[-2mm]
 Instituto Superior T\'{e}cnico,\\[-2mm]
Edif\'{\i}cio Ci\^{e}ncia,\\[-2mm] Piso 3,
1049-001 LISBOA, Portugal;\\[-2mm]
e-mail: isabel@cartan.ist.utl.pt}
\end{document}